\documentclass[a4paper]{amsart}
\usepackage[T1]{fontenc}
\usepackage[utf8]{inputenc}
\usepackage[active]{srcltx}
\usepackage{units}
\usepackage{amsthm}
\usepackage{amstext}
\usepackage{amssymb}
\usepackage{setspace}
\usepackage{esint}
\usepackage{xargs}[2008/03/08]

\makeatletter

\numberwithin{equation}{section}
\numberwithin{figure}{section}
\theoremstyle{plain}
\newtheorem{thm}{\protect\theoremname}
  \theoremstyle{plain}
  \newtheorem{prop}[thm]{\protect\propositionname}
  \theoremstyle{remark}
  \newtheorem*{rem*}{\protect\remarkname}

\makeatother

  \providecommand{\propositionname}{Proposition}
  \providecommand{\remarkname}{Remark}
\providecommand{\theoremname}{Theorem}

\begin{document}

\title{Hermite and poly-Bernoulli mixed-type Polynomials}

\author{Dae San Kim and Taekyun Kim}
\begin{abstract}
In this paper, we consider Hermite and poly-Bernoulli mixed-type polynomials
and investigate the properties of those polynomials which are derived
from umbral calculus. Finally, we give various identities associated
with Stirling numbers, Bernoulli and Frobenius-Euler polynomials of
higher order.

\newcommandx\li[1][usedefault, addprefix=\global, 1=]{\textnormal{Li}_{#1}}

\end{abstract}
\maketitle

\section{Introduction}

For $r\in\mathbb{Z}_{\ge0}$, as is well known, the Bernoulli polynomials
of order $r$ are defined by the generating function to be
\begin{equation}
\sum_{n=0}^{\infty}\frac{\mathbb{B}_{n}^{\left(r\right)}\left(x\right)}{n!}t^{n}=\left(\frac{t}{e^{t}-1}\right)^{r}e^{xt},\quad\left(\textrm{see [1-16]}\right).\label{eq:1}
\end{equation}

For $k\in\mathbb{Z}$, the polylogarithm is defined by
\begin{equation}
\li[k]\left(x\right)=\sum_{n=1}^{\infty}\frac{x^{n}}{n^{k}}.\label{eq:2}
\end{equation}

Note that $\li[1]\left(x\right)=-\log\left(1-x\right).$

The poly-Bernoulli polynomials are defined by the generating function
to be
\begin{equation}
\frac{\li[k]\left(1-e^{-t}\right)}{1-e^{-t}}e^{xt}=\sum_{n=0}^{\infty}B_{n}^{\left(k\right)}\left(x\right)\frac{t^{n}}{n!},\quad\left(\textrm{see }\left[5,8\right]\right).\label{eq:3}
\end{equation}

When $x=0$, $B_{n}^{\left(k\right)}=B_{n}^{\left(k\right)}\left(0\right)$
are called the poly-Bernoulli numbers (of index $k$).

For $\nu\left(\ne0\right)\in\mathbb{R}$, the Hermite polynomials of
order $\nu$ are given by the generating function to be
\begin{equation}
e^{-\frac{\nu t^{2}}{2}}e^{xt}=\sum_{n=0}^{\infty}H_{n}^{\left(\nu\right)}\left(x\right)\frac{t^{n}}{n!}\quad\left(\textrm{see \ensuremath{\left[6,12,13\right]}}\right).\label{eq:4}
\end{equation}

When $x=0$, $H_{n}^{\left(\nu\right)}=H_{n}^{\left(\nu\right)}\left(0\right)$
are called the Hermite numbers of order $\nu$.

In this paper, we consider the Hermite and poly-Bernoulli mixed-type
polynomials $HB_{n}^{\left(\nu,k\right)}\left(x\right)$ which are
defined by the generating function to be
\begin{equation}
e^{-\frac{\nu t^{2}}{2}}\frac{\li[k]\left(1-e^{-t}\right)}{1-e^{-t}}e^{xt}=\sum_{n=0}^{\infty}HB_{n}^{\left(\nu,k\right)}\left(x\right)\frac{t^{n}}{n!},\label{eq:5}
\end{equation}
where $k\in\mathbb{Z}$ and $\nu\left(\ne0\right)\in\mathbb{R}$.

When $x=0$, $HB_{n}^{\left(\nu,k\right)}=HB_{n}^{\left(\nu,k\right)}\left(0\right)$
are called the Hermite and poly-Bernoulli mixed-type numbers.

Let $\mathcal{F}$ be the set of all formal power series in the variable
$t$ over $\mathbb{C}$ as follows:
\begin{equation}
\mathcal{F}=\left\{ \left.f\left(t\right)=\sum_{k=0}^{\infty}a_{k}\frac{t^{k}}{k!}\right|a_{k}\in\mathbb{C}\right\} .\label{eq:6}
\end{equation}

Let $\mathbb{P}=\mathbb{C}\left[x\right]$ and $\mathbb{P}^{*}$ denote
the vector space of all linear functionals on $\mathbb{P}$.

$\left\langle \left.L\right|p\left(x\right)\right\rangle $ denotes
the action of the linear functional $L$ on the polynomial $p\left(x\right)$,
and we recall that the vector space operations on $\mathbb{P}^{*}$
are defined by $\left\langle \left.L+M\right|p\left(x\right)\right\rangle =\left\langle \left.L\right|p\left(x\right)\right\rangle +\left\langle \left.M\right|p\left(x\right)\right\rangle $,
$\left\langle \left.cL\right|p\left(x\right)\right\rangle =c\left\langle \left.L\right|p\left(x\right)\right\rangle $,
where $c$ is complex constant in $\mathbb{C}$. For $f\left(t\right)\in\mathcal{F}$,
let us define the linear functional on $\mathbb{P}$ by setting
\begin{equation}
\left\langle \left.f\left(t\right)\right|x^{n}\right\rangle =a_{n},\quad\left(n\ge0\right).\label{eq:7}
\end{equation}

Then, by (\ref{eq:6}) and (\ref{eq:7}), we get
\begin{equation}
\left\langle \left.t^{k}\right|x^{n}\right\rangle =n!\delta_{n,k},\quad\left(n,k\ge0\right),\label{eq:8}
\end{equation}
where $\delta_{n,k}$ is the Kronecker's symbol.

For $f_{L}\left(t\right)={\displaystyle \sum_{k=0}^{\infty}\frac{\left\langle \left.L\right|x^{k}\right\rangle }{k!}t^{k},}$
we have $\left\langle \left.f_{L}\left(t\right)\right|x^{n}\right\rangle =\left\langle \left.L\right|x^{n}\right\rangle .$
That is, $L=f_{L}\left(t\right).$ The map $L\longmapsto f_{L}\left(t\right)$
is a vector space isomorphism from $\mathbb{P}^{*}$ onto $\mathcal{F}$.
Henceforth, $\mathcal{F}$ denotes both the algebra of formal power
series in $t$ and the vector space of all linear functionals on $\mathbb{P}$,
and so an element $f\left(t\right)$ of $\mathcal{F}$ will be thought
of as both a formal power series and a linear functional. We call
$\mathcal{F}$ the umbral algebra and the umbral calculus is the study
of umbral algebra. The order $O\left(f\right)$ of the power series
$f\left(t\right)\ne0$ is the smallest integer for which $a_{k}$
does not vanish. If $O\left(f\right)=0$, then $f\left(t\right)$
is called an invertible series. If $O\left(f\right)=1$, then $f\left(t\right)$
is called a delta series. For $f\left(t\right),g\left(t\right)\in\mathcal{F}$,
we have
\begin{equation}
\left\langle \left.f\left(t\right)g\left(t\right)\right|p\left(x\right)\right\rangle =\left\langle \left.f\left(t\right)\right|g\left(t\right)p\left(x\right)\right\rangle =\left\langle \left.g\left(t\right)\right|f\left(t\right)p\left(x\right)\right\rangle .\label{eq:9}
\end{equation}

Let $f\left(t\right)\in\mathcal{F}$ and $p\left(x\right)\in\mathbb{P}$.
Then we have
\begin{equation}
f\left(t\right)=\sum_{k=0}^{\infty}\frac{\left\langle \left.f\left(t\right)\right|x^{k}\right\rangle }{k!}t^{k},\quad p\left(x\right)=\sum_{k=0}^{\infty}\frac{\left\langle \left.t^{k}\right|p\left(x\right)\right\rangle }{k!}x^{k},\quad\left(\textrm{see }\left[8,9,11,13,14\right]\right).\label{eq:10}
\end{equation}

By (\ref{eq:10}), we get
\begin{equation}
p^{\left(k\right)}\left(0\right)=\left\langle \left.t^{k}\right|p\left(x\right)\right\rangle =\left\langle \left.1\right|p^{\left(k\right)}\left(x\right)\right\rangle ,\label{eq:11}
\end{equation}
where $p^{\left(k\right)}\left(0\right)=\left.\frac{d^{k}p\left(x\right)}{dx^{k}}\right|_{x=0}.$

From (\ref{eq:11}), we have
\begin{equation}
t^{k}p\left(x\right)=p^{\left(k\right)}\left(x\right)=\frac{d^{k}p\left(x\right)}{dx^{k}},\quad\left(\textrm{see \ensuremath{\left[8,9,13\right]}}\right).\label{eq:12}
\end{equation}

By (\ref{eq:12}), we easily get
\begin{equation}
e^{yt}p\left(x\right)=p\left(x+y\right),\quad\left\langle \left.e^{yt}\right|p\left(x\right)\right\rangle =p\left(y\right).\label{eq:13}
\end{equation}

For $O\left(f\left(t\right)\right)=1$, $O\left(g\left(t\right)\right)=0$,
there exists a unique sequence $s_{n}\left(x\right)$ of polynomials
such that $\left\langle \left.g\left(t\right)f\left(t\right)^{k}\right|x^{n}\right\rangle =n!\delta_{n,k},\quad\left(n,k\ge0\right)$.

The sequence $s_{n}\left(x\right)$ is called the Sheffer sequence for
$\left(g\left(t\right),f\left(t\right)\right)$ which is denoted by
$s_{n}\left(x\right)\sim\left(g\left(t\right),f\left(t\right)\right).$

Let $p\left(x\right)\in\mathbb{P},$ $f\left(t\right)\in\mathcal{F}$.
Then we see that
\begin{equation}
\left\langle \left.f\left(t\right)\right|xp\left(x\right)\right\rangle =\left\langle \left.\partial_{t}f(t)\right|p\left(x\right)\right\rangle =\left\langle \left.\frac{df\left(t\right)}{dt}\right|p\left(x\right)\right\rangle .\label{eq:14}
\end{equation}

For $s_{n}\left(x\right)\sim\left(g\left(t\right),f\left(t\right)\right),$
we have the following equations :
\begin{equation}
h\left(t\right)=\sum_{k=0}^{\infty}\frac{\left\langle \left.h\left(t\right)\right|s_{k}\left(x\right)\right\rangle }{k!}g\left(t\right)f\left(t\right)^{k},\quad p\left(x\right)=\sum_{k=0}^{\infty}\frac{\left\langle \left.g\left(t\right)f\left(t\right)^{k}\right|p\left(x\right)\right\rangle }{k!}s_{k}\left(x\right),\label{eq:15}
\end{equation}
where $h\left(t\right)\in\mathcal{F},$ $p\left(x\right)\in\mathbb{P}$.

\begin{equation}
\frac{1}{g\left(\overline{f}\left(t\right)\right)}e^{y\overline{f}\left(t\right)}=\sum_{n=0}^{\infty}s_{n}\left(y\right)\frac{t^{n}}{n!},\label{eq:16}
\end{equation}
where $\overline{f}\left(t\right)$ is the compositional inverse for
$f\left(t\right)$ with $f\left(\overline{f}\left(t\right)\right)=t$,
\begin{equation}
s_{n}\left(x+y\right)=\sum_{k=0}^{n}\dbinom{n}{k}s_{k}\left(y\right)p_{n-k}\left(x\right),\quad\textrm{where }p_{n}\left(x\right)=g\left(t\right)s_{n}\left(x\right),\label{eq:17}
\end{equation}
\begin{equation}
f\left(t\right)s_{n}\left(x\right)=ns_{n-1}\left(x\right),\quad s_{n+1}\left(x\right)=\left(x-\frac{g^{\prime}\left(t\right)}{g\left(t\right)}\right)\frac{1}{f^{\prime}(t)}s_{n}\left(x\right),\label{eq:18}
\end{equation}
 and the conjugate representation is given by
\begin{equation}
s_{n}\left(x\right)=\sum_{j=0}^{n}\frac{1}{j!}\left\langle \left.g\left(\overline{f}\left(t\right)\right)^{-1}\overline{f}\left(t\right)^{j}|x^{n}\right.\right\rangle x^{j}.\label{eq:19}
\end{equation}

For $s_{n}\left(x\right)\sim\left(g\left(t\right),f\left(t\right)\right)$,
$r_{n}\left(x\right)\sim\left(h\left(t\right),l\left(t\right)\right)$,
we have
\begin{equation}
s_{n}\left(x\right)=\sum_{m=0}^{n}C_{n,m}r_{m}\left(x\right),\label{eq:20}
\end{equation}
 where
\begin{equation}
C_{n,m}=\frac{1}{m!}\left\langle \left.\frac{h\left(\overline{f}\left(t\right)\right)}{g\left(\overline{f}\left(t\right)\right)}l\left(\overline{f}\left(t\right)\right)^{m}\right|x^{n}\right\rangle ,\quad\left(\textrm{see \ensuremath{\left[8,9,13\right]}}\right).\label{eq:21}
\end{equation}

In this paper, we consider Hermite and poly-Bernoulli mixed-type polynomials
and investigate the properties of those polynomials which are derived
from umbral calculus. Finally, we give various identities associated
with Bernoulli and Frobenius-Euler polynomials of higher order.

\section{Hermite and poly-Bernoulli mixed-type polynomials}

$\,$

From (\ref{eq:5}) and (\ref{eq:16}), we note that
\begin{equation}
HB_{n}^{\left(\nu,k\right)}\left(x\right)\sim\left(e^{\frac{\nu t^{2}}{2}}\frac{1-e^{-t}}{\li[k]\left(1-e^{-t}\right)},t\right),\label{eq:22}
\end{equation}
 and, by (\ref{eq:3}), (\ref{eq:4}) and (\ref{eq:16}), we get
\begin{equation}
B_{n}^{\left(k\right)}\left(x\right)\sim\left(\frac{1-e^{-t}}{\li[k]\left(1-e^{-t}\right)},t\right),\label{eq:23}
\end{equation}

\begin{equation}
H_{n}^{\left(\nu\right)}\left(x\right)\sim\left(e^{\frac{\nu t^{2}}{2}},t\right),\quad\textrm{where }n\ge0.\label{eq:24}
\end{equation}

From (\ref{eq:18}), (\ref{eq:22}), (\ref{eq:23}), and (\ref{eq:24}),
we have
\begin{equation}
tB_{n}^{\left(k\right)}\left(x\right)=nB_{n-1}^{\left(k\right)}\left(x\right),\quad tH_{n}^{\left(\nu\right)}\left(x\right)=nH_{n-1}^{\left(\nu\right)}\left(x\right),\quad tHB_{n}^{\left(\nu,k\right)}\left(x\right)=nHB_{n-1}^{\left(\nu,k\right)}\left(x\right).\label{eq:25}
\end{equation}

By (\ref{eq:5}), (\ref{eq:8}) and (\ref{eq:22}), we get

\begin{eqnarray}
HB_{n}^{\left(\nu,k\right)}\left(x\right) & = & e^{-\frac{\nu t^{2}}{2}}\frac{\li[k]\left(1-e^{-t}\right)}{1-e^{-t}}x^{n}=e^{-\frac{\nu t^{2}}{2}}B_{n}^{\left(k\right)}\left(x\right)\label{eq:26}\\
 & = & \sum_{m=0}^{\left[\frac{n}{2}\right]}\frac{1}{m!}\left(-\frac{\nu}{2}\right)^{m}\left(n\right)_{2m}B_{n-2m}^{\left(k\right)}\left(x\right)\nonumber \\
 & = & \sum_{m=0}^{\left[\frac{n}{2}\right]}\dbinom{n}{2m}\frac{\left(2m\right)!}{m!}\left(-\frac{\nu}{2}\right)^{m}B_{n-2m}^{\left(k\right)}\left(x\right).\nonumber
\end{eqnarray}

Therefore, by (\ref{eq:26}), we obtain the following proposition.

\begin{prop}

For $n\ge0$, we have
\[
HB_{n}^{\left(\nu,k\right)}\left(x\right)=\sum_{m=0}^{\left[\frac{n}{2}\right]}\dbinom{n}{2m}\frac{\left(2m\right)!}{m!}\left(-\frac{\nu}{2}\right)^{m}B_{n-2m}^{\left(k\right)}\left(x\right).
\]

\end{prop}

$\,$

From (\ref{eq:5}), we can also derive
\begin{eqnarray}
HB_{n}^{\left(\nu,k\right)}\left(x\right) & = & \frac{\li[k]\left(1-e^{-t}\right)}{1-e^{-t}}e^{-\frac{\nu t^{2}}{2}}x^{n}=\frac{\li[k]\left(1-e^{-t}\right)}{1-e^{-t}}H_{n}^{\left(\nu\right)}\left(x\right)\label{eq:27}\\
 & = & \sum_{m=0}^{\infty}\frac{\left(1-e^{-t}\right)^{m}}{\left(m+1\right)^{k}}H_{n}^{\left(\nu\right)}\left(x\right)\nonumber \\
 & = & \sum_{m=0}^{n}\frac{1}{\left(m+1\right)^{k}}\sum_{j=0}^{m}\dbinom{m}{j}\left(-1\right)^{j}e^{-jt}H_{n}^{\left(\nu\right)}\left(x\right)\nonumber \\
 & = & \sum_{m=0}^{n}\frac{1}{\left(m+1\right)^{k}}\sum_{j=0}^{m}\dbinom{m}{j}\left(-1\right)^{j}H_{n}^{\left(\nu\right)}\left(x-j\right).\nonumber
\end{eqnarray}

Therefore, by (\ref{eq:27}), we obtain the following theorem.

\begin{thm}

For $n\ge0$, we have
\[
HB_{n}^{\left(\nu,k\right)}\left(x\right)=\sum_{m=0}^{n}\frac{1}{\left(m+1\right)^{k}}\sum_{j=0}^{m}\dbinom{m}{j}\left(-1\right)^{j}H_{n}^{\left(\nu\right)}\left(x-j\right).
\]

\end{thm}

$\,$

By (\ref{eq:5}), we get

\begin{align}
 & HB_{n}^{\left(\nu,k\right)}\left(x\right)\label{eq:28}\\
= & e^{-\frac{\nu t^{2}}{2}}B_{n}^{\left(k\right)}\left(x\right)=\sum_{l=0}^{\infty}\frac{1}{l!}\left(-\frac{\nu}{2}\right)^{l}t^{2l}B_{n}^{\left(k\right)}\left(x\right)\nonumber \\
= & \sum_{l=0}^{\left[\frac{n}{2}\right]}\frac{1}{l!}\left(-\frac{\nu}{2}\right)^{l}\sum_{m=0}^{n}\frac{1}{\left(m+1\right)^{k}}\sum_{j=0}^{m}\left(-1\right)^{j}\dbinom{m}{j}t^{2l}\left(x-j\right)^{n}\nonumber \\
= & \sum_{l=0}^{\left[\frac{n}{2}\right]}\sum_{j=0}^{n}\left\{ \sum_{m=j}^{n}\dbinom{n}{2l}\frac{\left(2l\right)!}{l!}\left(-\frac{\nu}{2}\right)^{l}\frac{\left(-1\right)^{j}\tbinom{m}{j}}{\left(m+1\right)^{k}}\right\} \left(x-j\right)^{n-2l}.\nonumber
\end{align}

Therefore, by (\ref{eq:28}), we obtain the following theorem.

\begin{thm}

For $n\ge0$, we have
\[
HB_{n}^{\left(\nu,k\right)}\left(x\right)=\sum_{l=0}^{\left[\frac{n}{2}\right]}\sum_{j=0}^{n}\left\{ \sum_{m=j}^{n}\dbinom{n}{2l}\frac{\left(2l\right)!}{l!}\left(-\frac{\nu}{2}\right)^{l}\frac{\left(-1\right)^{j}\tbinom{m}{j}}{\left(m+1\right)^{k}}\right\} \left(x-j\right)^{n-2l}.
\]

\end{thm}

$\,$

By (\ref{eq:27}), we get
\begin{align}
 & HB_{n}^{\left(\nu,k\right)}\left(x\right)\label{eq:29}\\
= & \sum_{m=0}^{n}\frac{\left(1-e^{-t}\right)^{m}}{\left(m+1\right)^{k}}H_{n}^{\left(\nu\right)}\left(x\right)\nonumber \\
= & \sum_{m=0}^{n}\frac{1}{\left(m+1\right)^{k}}\sum_{a=0}^{n-m}\frac{m!}{\left(a+m\right)!}\left(-1\right)^{a}S_{2}\left(a+m,m\right)\left(n\right)_{a+m}H_{n-a-m}^{\left(\nu\right)}\left(x\right)\nonumber \\
= & \sum_{m=0}^{n}\sum_{a=0}^{n-m}\frac{\left(-1\right)^{n-a-m}m!}{\left(m+1\right)^{k}}\dbinom{n}{n-a}S_{2}\left(n-a,m\right)H_{a}^{\left(\nu\right)}\left(x\right)\nonumber \\
= & \left(-1\right)^{n}\sum_{a=0}^{n}\left\{ \sum_{m=0}^{n-a}\frac{\left(-1\right)^{m+a}m!}{\left(m+1\right)^{k}}\dbinom{n}{a}S_{2}\left(n-a,m\right)\right\} H_{a}^{\left(\nu\right)}\left(x\right),\nonumber
\end{align}
 where $S_{2}\left(n,m\right)$ is the Stirling number of the second
kind.

Therefore, by (\ref{eq:29}), we obtain the following theorem.

\begin{thm}

For $n\ge0$, we have
\[
HB_{n}^{\left(\nu,k\right)}\left(x\right)=\left(-1\right)^{n}\sum_{a=0}^{n}\left\{ \sum_{m=0}^{n-a}\frac{\left(-1\right)^{a+m}m!}{\left(m+1\right)^{k}}\dbinom{n}{a}S_{2}\left(n-a,m\right)\right\} H_{a}^{\left(\nu\right)}\left(x\right).
\]

\end{thm}

$\,$

From (\ref{eq:19}) and (\ref{eq:22}), we have
\begin{eqnarray}
HB_{n}^{\left(\nu,k\right)}\left(x\right) & = & \sum_{j=0}^{n}\dbinom{n}{j}\left\langle \left.e^{-\frac{\nu t^{2}}{2}}\frac{\li[k]\left(1-e^{-t}\right)}{1-e^{-t}}\right|x^{n-j}\right\rangle x^{j}\label{eq:30}\\
 & = & \sum_{j=0}^{n}\dbinom{n}{j}\left\langle \left.e^{-\frac{\nu t^{2}}{2}}\right|B_{n-j}^{\left(k\right)}\left(x\right)\right\rangle x^{j}\nonumber \\
 & = & \sum_{j=0}^{n}\dbinom{n}{j}\sum_{l=0}^{\left[\frac{n-j}{2}\right]}\frac{\left(-\frac{\nu}{2}\right)^{l}}{l!}\left(n-j\right)_{2l}\left\langle \left.1\right|B_{n-j-2l}^{\left(k\right)}\left(x\right)\right\rangle x^{j}\nonumber \\
 & = & \sum_{j=0}^{n}\dbinom{n}{j}\sum_{l=0}^{\left[\frac{n-j}{2}\right]}\frac{1}{l!}\left(-\frac{\nu}{2}\right)^{l}\left(n-j\right)_{2l}B_{n-j-2l}^{\left(k\right)}x^{j}\nonumber \\
 & = & \sum_{j=0}^{n}\left\{ \sum_{l=0}^{\left[\frac{n-j}{2}\right]}\dbinom{n}{j}\dbinom{n-j}{2l}\frac{\left(2l\right)!}{l!}\left(-\frac{\nu}{2}\right)^{l}B_{n-j-2l}^{\left(k\right)}\right\} x^{j}.\nonumber
\end{eqnarray}

Therefore, by (\ref{eq:30}), we obtain the following theorem.

\begin{thm}

For $n\ge0$, we have
\[
HB_{n}^{\left(\nu,k\right)}\left(x\right)=\sum_{j=0}^{n}\left\{ \sum_{l=0}^{\left[\frac{n-j}{2}\right]}\dbinom{n}{j}\dbinom{n-j}{2l}\frac{\left(2l\right)!}{l!}\left(-\frac{\nu}{2}\right)^{l}B_{n-j-2l}^{\left(k\right)}\right\} x^{j}.
\]

\end{thm}
\begin{rem*}

By (\ref{eq:17}) and (\ref{eq:22}), we easily get
\begin{equation}
HB_{n}^{\left(\nu,k\right)}\left(x+y\right)=\sum_{j=0}^{n}\dbinom{n}{j}HB_{j}^{\left(\nu,k\right)}\left(x\right)y^{n-j}.\label{eq:31}
\end{equation}

\end{rem*}

$\,$

We note that
\begin{equation}
HB_{n}^{\left(\nu,k\right)}\left(x\right)\sim\left(g\left(t\right)=e^{\frac{\nu t^{2}}{2}}\frac{1-e^{-t}}{\li[k]\left(1-e^{-t}\right)},\, f\left(t\right)=t\right).\label{eq:32}
\end{equation}

From (\ref{eq:18}) and (\ref{eq:32}), we have

\begin{equation}
HB_{n+1}^{\left(\nu,k\right)}\left(x\right)=\left(x-\frac{g^{\prime}\left(t\right)}{g\left(t\right)}\right)HB_{n}^{\left(\nu,k\right)}\left(x\right).\label{eq:33}
\end{equation}

Now, we observe that
\begin{eqnarray}
\frac{g^{\prime}\left(t\right)}{g\left(t\right)} & = & \left(\log\left(g\left(t\right)\right)\right)^{\prime}\label{eq:34}\\
 & = & \left(\log e^{\frac{\nu t^{2}}{2}}+\log\left(1-e^{-t}\right)-\log\left(\li[k]\left(1-e^{-t}\right)\right)\right)^{\prime}\nonumber \\
 & = & \nu t+\frac{e^{-t}}{1-e^{-t}}\left(1-\frac{\li[k-1]\left(1-e^{-t}\right)}{\li[k]\left(1-e^{-t}\right)}\right).\nonumber
\end{eqnarray}

By (\ref{eq:33}) and (\ref{eq:34}), we get
\begin{align}
 & HB_{n+1}^{\left(\nu,k\right)}\left(x\right)\label{eq:35}\\
= & xHB_{n}^{\left(\nu,k\right)}\left(x\right)-\frac{g^{\prime}\left(t\right)}{g\left(t\right)}HB_{n}^{\left(\nu,k\right)}\left(x\right)\nonumber \\
= & xHB_{n}^{\left(\nu,k\right)}\left(x\right)-\nu nHB_{n-1}^{\left(\nu,k\right)}\left(x\right)-e^{-\frac{\nu t^{2}}{2}}\frac{t}{e^{t}-1}\frac{\li[k]\left(1-e^{-t}\right)-\li[k-1]\left(1-e^{-t}\right)}{t\left(1-e^{-t}\right)}x^{n}.\nonumber
\end{align}

It is easy to show that
\begin{align}
 & \frac{\li[k]\left(1-e^{-t}\right)-\li[k-1]\left(1-e^{-t}\right)}{1-e^{-t}}\label{eq:36}\\
= & \sum_{m=2}^{\infty}\left(\frac{1}{m^{k}}-\frac{1}{m^{k-1}}\right)\left(1-e^{-t}\right)^{m-1}\nonumber \\
= & \left(\frac{1}{2^{k}}-\frac{1}{2^{k-1}}\right)t+\cdots.\nonumber
\end{align}

Thus, by (\ref{eq:36}), we get
\begin{equation}
\frac{\li[k]\left(1-e^{-t}\right)-\li[k-1]\left(1-e^{-t}\right)}{t\left(1-e^{-t}\right)}x^{n}=\frac{\li[k]\left(1-e^{-t}\right)-\li[k-1]\left(1-e^{-t}\right)}{1-e^{-t}}\frac{x^{n+1}}{n+1}.\label{eq:37}
\end{equation}

From (\ref{eq:37}), we can derive
\begin{align}
 & e^{-\frac{\nu t^{2}}{2}}\frac{t}{e^{t}-1}\frac{\li[k]\left(1-e^{-t}\right)-\li[k-1]\left(1-e^{-t}\right)}{t\left(1-e^{-t}\right)}x^{n}\label{eq:38}\\
= & \frac{1}{n+1}\left(\sum_{l=0}^{\infty}\frac{B_{l}}{l!}t^{l}\right)\left(HB_{n+1}^{\left(\nu,k\right)}\left(x\right)-HB_{n+1}^{\left(\nu,k-1\right)}\left(x\right)\right)\nonumber \\
= & \frac{1}{n+1}\sum_{l=0}^{n+1}\frac{B_{l}}{l!}t^{l}\left(HB_{n+1}^{\left(\nu,k\right)}\left(x\right)-HB_{n+1}^{\left(\nu,k-1\right)}\left(x\right)\right)\nonumber \\
= & \frac{1}{n+1}\sum_{l=0}^{n+1}\dbinom{n+1}{l}B_{l}\left(HB_{n+1-l}^{\left(\nu,k\right)}\left(x\right)-HB_{n+1-l}^{\left(\nu,k-1\right)}\left(x\right)\right).\nonumber
\end{align}

Therefore, by (\ref{eq:35}) and (\ref{eq:38}), we obtain the following
theorem.

\begin{thm}

For $n\ge0$, we have
\begin{align}
 & HB_{n+1}^{\left(\nu,k\right)}\left(x\right)\label{eq:39}\\
= & xHB_{n}^{\left(\nu,k\right)}\left(x\right)-\nu nHB_{n-1}^{\left(\nu,k\right)}\left(x\right)\nonumber \\
 & -\frac{1}{n+1}\sum_{l=0}^{n+1}\dbinom{n+1}{l}B_{l}\left\{ HB_{n+1-l}^{\left(\nu,k\right)}\left(x\right)-HB_{n+1-l}^{\left(\nu,k-1\right)}\left(x\right)\right\} .\nonumber
\end{align}

\end{thm}

Let us take $t$ on the both sides of (\ref{eq:39}). Then we have
\begin{align}
 & \left(n+1\right)HB_{n}^{\left(\nu,k\right)}\left(x\right)\label{eq:40}\\
= & \left(xt+1\right)HB_{n}^{\left(\nu,k\right)}\left(x\right)-\nu n\left(n-1\right)HB_{n-2}^{\left(\nu,k\right)}\left(x\right)\nonumber \\
 & -\frac{1}{n+1}\sum_{l=0}^{n+1}\dbinom{n+1}{l}\left(n+1-l\right)B_{l}\left\{ HB_{n-l}^{\left(\nu,k\right)}\left(x\right)-HB_{n-l}^{\left(\nu,k-1\right)}\left(x\right)\right\} \nonumber \\
= & nxHB_{n-1}^{\left(\nu,k\right)}\left(x\right)+HB_{n}^{\left(\nu,k\right)}\left(x\right)-\nu n\left(n-1\right)HB_{n-2}^{\left(\nu,k\right)}\left(x\right)\nonumber \\
 & -\sum_{l=0}^{n}\dbinom{n}{l}B_{l}\left(HB_{n-l}^{\left(\nu,k\right)}\left(x\right)-HB_{n-l}^{\left(\nu,k-1\right)}\left(x\right)\right),\nonumber
\end{align}
where $n\ge3$.

Thus, by (\ref{eq:40}), we obtain the following theorem.

\begin{thm}

For $n\ge3$, we have
\begin{align*}
 & \sum_{l=0}^{n}\dbinom{n}{l}B_{l}HB_{n-l}^{\left(\nu,k-1\right)}\left(x\right)\\
= & \left(n+1\right)HB_{n}^{\left(\nu,k\right)}\left(x\right)-n\left(x+\frac{1}{2}\right)HB_{n-1}^{\left(\nu,k\right)}\left(x\right)\\
 & +n\left(n-1\right)\left(\nu+\frac{1}{12}\right)HB_{n-2}^{\left(\nu,k\right)}\left(x\right)\\
 & +\sum_{l=0}^{n-3}\dbinom{n}{l}B_{n-l}HB_{l}^{\left(\nu,k\right)}\left(x\right).
\end{align*}

\end{thm}

By (\ref{eq:5}) and (\ref{eq:8}), we get

\begin{align}
 & HB_{n}^{\left(\nu,k\right)}\left(y\right)\label{eq:41}\\
= & \left\langle \left.e^{-\frac{\nu t^{2}}{2}}\frac{\li[k]\left(1-e^{-t}\right)}{1-e^{-t}}e^{yt}\right|x^{n}\right\rangle \nonumber \\
= & \left\langle \left.\partial_{t}\left(e^{-\frac{\nu t^{2}}{2}}\frac{\li[k]\left(1-e^{-t}\right)}{1-e^{-t}}e^{yt}\right)\right|x^{n-1}\right\rangle \nonumber \\
= & \left\langle \left.\left(\partial_{t}e^{-\frac{\nu t^{2}}{2}}\right)\frac{\li[k]\left(1-e^{-t}\right)}{1-e^{-t}}e^{yt}\right|x^{n-1}\right\rangle \nonumber \\
 & +\left\langle \left.e^{-\frac{\nu t^{2}}{2}}\left(\partial_{t}\frac{\li[k]\left(1-e^{-t}\right)}{1-e^{-t}}\right)e^{yt}\right|x^{n-1}\right\rangle \nonumber \\
 & +\left\langle \left.e^{-\frac{\nu t^{2}}{2}}\frac{\li[k]\left(1-e^{-t}\right)}{1-e^{-t}}\left(\partial_{t}e^{yt}\right)\right|x^{n-1}\right\rangle \nonumber
\end{align}

\begin{align*}
= & -\nu\left(n-1\right)\left\langle \left.e^{-\frac{\nu t^{2}}{2}}\frac{\li[k]\left(1-e^{-t}\right)}{1-e^{-t}}e^{yt}\right|x^{n-2}\right\rangle \\
 & +y\left\langle \left.e^{-\frac{\nu t^{2}}{2}}\frac{\li[k]\left(1-e^{-t}\right)}{1-e^{-t}}e^{yt}\right|x^{n-1}\right\rangle \\
 & +\left\langle \left.e^{-\frac{\nu t^{2}}{2}}\left(\partial_{t}\frac{\li[k]\left(1-e^{-t}\right)}{1-e^{-t}}\right)e^{yt}\right|x^{n-1}\right\rangle \\
= & -\nu\left(n-1\right)HB_{n-2}^{\left(\nu,k\right)}\left(y\right)+yHB_{n-1}^{\left(\nu,k\right)}\left(y\right)\\
 & +\left\langle \left.e^{-\frac{\nu t^{2}}{2}}\left(\partial_{t}\frac{\li[k]\left(1-e^{-t}\right)}{1-e^{-t}}\right)e^{yt}\right|x^{n-1}\right\rangle .
\end{align*}

Now, we observe that

\begin{equation}
\partial_{t}\left(\frac{\li[k]\left(1-e^{-t}\right)}{1-e^{-t}}\right)=\frac{\li[k-1]\left(1-e^{-t}\right)-\li[k]\left(1-e^{-t}\right)}{\left(1-e^{-t}\right)^{2}}e^{-t}.\label{eq:42}
\end{equation}

From (\ref{eq:42}), we have
\begin{align}
 & \left\langle \left.e^{-\frac{\nu t^{2}}{2}}\left(\partial_{t}\frac{\li[k]\left(1-e^{-t}\right)}{1-e^{-t}}\right)e^{yt}\right|x^{n-1}\right\rangle \label{eq:43}\\
= & \left\langle \left.e^{-\frac{\nu t^{2}}{2}}\left(\frac{\li[k-1]\left(1-e^{-t}\right)-\li[k]\left(1-e^{-t}\right)}{\left(1-e^{-t}\right)^{2}}\right)e^{-t}e^{yt}\right|\frac{1}{n}tx^{n}\right\rangle \nonumber \\
= & \frac{1}{n}\left\langle \left.e^{-\frac{\nu t^{2}}{2}}\frac{\li[k-1]\left(1-e^{-t}\right)-\li[k]\left(1-e^{-t}\right)}{1-e^{-t}}e^{yt}\right|\frac{t}{e^{t}-1}x^{n}\right\rangle \nonumber \\
= & \frac{1}{n}\left\langle \left.e^{-\frac{\nu t^{2}}{2}}\frac{\li[k-1]\left(1-e^{-t}\right)-\li[k]\left(1-e^{-t}\right)}{1-e^{-t}}e^{yt}\right|B_{n}\left(x\right)\right\rangle \nonumber \\
= & \frac{1}{n}\sum_{l=0}^{n}\dbinom{n}{l}B_{l}\left\langle \left.e^{-\frac{\nu t^{2}}{2}}\frac{\li[k-1]\left(1-e^{-t}\right)-\li[k]\left(1-e^{-t}\right)}{1-e^{-t}}e^{yt}\right|x^{n-l}\right\rangle \nonumber \\
= & \frac{1}{n}\sum_{l=0}^{n}\dbinom{n}{l}B_{l}\left\{ HB_{n-l}^{\left(\nu,k-1\right)}\left(y\right)-HB_{n-l}^{\left(\nu,k\right)}\left(y\right)\right\} ,\nonumber
\end{align}
where $B_{n}$ are the ordinary Bernoulli numbers which are
defined by the generating function to be
\[
\frac{t}{e^{t}-1}=\sum_{n=0}^{\infty}\frac{B_{n}}{n!}t^{n}.
\]

Therefore, by (\ref{eq:41}) and (\ref{eq:43}), we obtain the following
theorem.

\begin{thm}

For $n\ge2$, we have
\begin{eqnarray*}
HB_{n}^{\left(\nu,k\right)}\left(x\right) & = & -\nu\left(n-1\right)HB_{n-2}^{\left(\nu,k\right)}\left(x\right)+xHB_{n-1}^{\left(\nu,k\right)}\left(x\right)\\
 &  & +\frac{1}{n}\sum_{l=0}^{n}\dbinom{n}{l}B_{l}\left(HB_{n-l}^{\left(\nu,k-1\right)}\left(x\right)-HB_{n-l}^{\left(\nu,k\right)}\left(x\right)\right).
\end{eqnarray*}

\end{thm}

$\,$

Now, we compute
\[
\left\langle \left.e^{-\frac{\nu t^{2}}{2}}\li[k]\left(1-e^{-t}\right)\right|x^{n+1}\right\rangle
\]
in two different ways.

On the one hand,
\begin{align}
 & \left\langle \left.e^{-\frac{\nu t^{2}}{2}}\li[k]\left(1-e^{-t}\right)\right|x^{n+1}\right\rangle \label{eq:45}\\
= & \left\langle \left.e^{-\frac{\nu t^{2}}{2}}\frac{\li[k]\left(1-e^{-t}\right)}{1-e^{-t}}\left(1-e^{-t}\right)\right|x^{n+1}\right\rangle \nonumber \\
= & \left\langle \left.e^{-\frac{\nu t^{2}}{2}}\frac{\li[k]\left(1-e^{-t}\right)}{1-e^{-t}}\right|\left(1-e^{-t}\right)x^{n+1}\right\rangle \nonumber \\
= & \left\langle \left.e^{-\frac{\nu t^{2}}{2}}\frac{\li[k]\left(1-e^{-t}\right)}{1-e^{-t}}\right|x^{n+1}-\left(x-1\right)^{n+1}\right\rangle \nonumber \\
= & \sum_{m=0}^{n}\left(-1\right)^{n-m}\dbinom{n+1}{m}\left\langle \left.e^{-\frac{\nu t^{2}}{2}}\frac{\li[k]\left(1-e^{-t}\right)}{1-e^{-t}}\right|x^{m}\right\rangle \nonumber \\
= & \sum_{m=0}^{n}\left(-1\right)^{n-m}\dbinom{n+1}{m}HB_{m}^{\left(\nu,k\right)}.\nonumber
\end{align}

On the other hand,
\begin{align}
 & \left\langle \left.e^{-\frac{\nu t^{2}}{2}}\li[k]\left(1-e^{-t}\right)\right|x^{n+1}\right\rangle \label{eq:46}\\
= & \left\langle \left.\li[k]\left(1-e^{-t}\right)\right|e^{-\frac{\nu t^{2}}{2}}x^{n+1}\right\rangle \nonumber \\
= & \left\langle \left.\int_{0}^{t}\left(\li[k]\left(1-e^{-s}\right)\right)^{\prime}ds\right|e^{-\frac{\nu t^{2}}{2}}x^{n+1}\right\rangle \nonumber \\
= & \left\langle \left.\int_{0}^{t}e^{-s}\frac{\li[k-1]\left(1-e^{-s}\right)}{1-e^{-s}}ds\right|e^{-\frac{\nu t^{2}}{2}}x^{n+1}\right\rangle \nonumber \\
= & \left\langle \left.\sum_{l=0}^{\infty}\left(\sum_{m=0}^{l}\left(-1\right)^{l-m}\dbinom{l}{m}B_{m}^{\left(k-1\right)}\frac{t^{l+1}}{\left(l+1\right)!}\right)\right|H_{n+1}^{\left(\nu\right)}\left(x\right)\right\rangle \nonumber \\
= & \sum_{l=0}^{n}\sum_{m=0}^{l}\left(-1\right)^{l-m}\dbinom{l}{m}B_{m}^{\left(k-1\right)}\frac{1}{\left(l+1\right)!}\left\langle \left.t^{l+1}\right|H_{n+1}^{\left(\nu\right)}\left(x\right)\right\rangle \nonumber \\
= & \sum_{l=0}^{n}\sum_{m=0}^{l}\left(-1\right)^{l-m}\dbinom{l}{m}\dbinom{n+1}{l+1}B_{m}^{\left(k-1\right)}H_{n-l}^{\left(\nu\right)}.\nonumber
\end{align}

Therefore, by (\ref{eq:45}) and (\ref{eq:46}), we obtain the following
theorem.

\begin{thm}

For $n\ge0$, we have

\begin{align*}
 & \sum_{m=0}^{n}\left(-1\right)^{n-m}\dbinom{n+1}{m}HB_{m}^{\left(\nu,k\right)}\\
= & \sum_{m=0}^{n}\sum_{l=m}^{n}\left(-1\right)^{l-m}\dbinom{l}{m}\dbinom{n+1}{l+1}B_{m}^{\left(k-1\right)}H_{n-l}^{\left(\nu\right)}.
\end{align*}

\end{thm}

$\,$

Let us consider the following two Sheffer sequences:

\begin{equation}
HB_{n}^{\left(\nu,k\right)}\left(x\right)\sim\left(e^{\frac{\nu t^{2}}{2}}\frac{1-e^{-t}}{\li[k]\left(1-e^{-t}\right)},t\right),\label{eq:47}
\end{equation}
and
\begin{equation}
\mathbb{B}_{n}^{\left(r\right)}\left(x\right)\sim\left(\left(\frac{e^{t}-1}{t}\right)^{r},t\right),\quad\left(r\in\mathbb{Z}_{\ge0}\right).\label{eq:48}
\end{equation}

Let us assume that
\begin{equation}
HB_{n}^{\left(\nu,k\right)}\left(x\right)=\sum_{m=0}^{n}C_{n,m}\mathbb{B}_{m}^{\left(r\right)}\left(x\right).\label{eq:49}
\end{equation}

Then, by (\ref{eq:20}) and (\ref{eq:21}), we get
\begin{eqnarray}
C_{n,m} & = & \frac{1}{m!}\left\langle \left.\left(\frac{e^{t}-1}{t}\right)^{r}t^{m}\right|e^{-\frac{\nu t^{2}}{2}}\frac{\li[k]\left(1-e^{-t}\right)}{1-e^{-t}}x^{n}\right\rangle \label{eq:50}\\
 & = & \frac{1}{m!}\left\langle \left.\left(\frac{e^{t}-1}{t}\right)^{r}\right|t^{m}HB_{n}^{\left(\nu,k\right)}\left(x\right)\right\rangle \nonumber \\
 & = & \frac{1}{m!}\left(n\right)_{m}\left\langle \left.\left(\frac{e^{t}-1}{t}\right)^{r}\right|HB_{n-m}^{\left(\nu,k\right)}\left(x\right)\right\rangle \nonumber \\
 & = & \dbinom{n}{m}\sum_{l=0}^{\infty}\frac{r!}{\left(l+r\right)!}S_{2}\left(l+r,r\right)\left\langle \left.t^{l}\right|HB_{n-m}^{\left(\nu,k\right)}\left(x\right)\right\rangle \nonumber \\
 & = & \dbinom{n}{m}\sum_{l=0}^{n-m}\left(n-m\right)_{l}\frac{r!}{\left(l+r\right)!}S_{2}\left(l+r,r\right)HB_{n-m-l}^{\left(\nu,k\right)}\nonumber \\
 & = & \dbinom{n}{m}\sum_{l=0}^{n-m}\frac{\tbinom{n-m}{l}}{\tbinom{l+r}{r}}S_{2}\left(l+r,r\right)HB_{n-m-l}^{\left(\nu,k\right)}.\nonumber
\end{eqnarray}

Therefore, by (\ref{eq:49}) and (\ref{eq:50}), we obtain the following
theorem.

\begin{thm}

For $n,r\in\mathbb{Z}_{\ge0}$, we have
\[
HB_{n}^{\left(\nu,k\right)}\left(x\right)=\sum_{m=0}^{n}\left\{ \dbinom{n}{m}\sum_{l=0}^{n-m}\frac{\tbinom{n-m}{l}}{\tbinom{l+r}{r}}S_{2}\left(l+r,r\right)HB_{n-m-l}^{\left(\nu,k\right)}\right\} \mathbb{B}_{m}^{\left(r\right)}\left(x\right).
\]

\end{thm}
$\,$

For $\lambda\left(\ne1\right)\in\mathbb{C}$, $r\in\mathbb{Z}_{\ge0}$,
the Frobenius-Euler polynomials of order $r$ are defined by the generating
function to be
\begin{equation}
\left(\frac{1-\lambda}{e^{t}-\lambda}\right)^{r}e^{xt}=\sum_{n=0}^{\infty}H_{n}^{\left(r\right)}\left(x|\lambda\right)\frac{t^{n}}{n!},\quad\left(\textrm{see \ensuremath{\left[1,4,7,9,10\right]}}\right).\label{eq:51}
\end{equation}

From (\ref{eq:16}) and (\ref{eq:51}), we note that
\begin{equation}
H_{n}^{\left(r\right)}\left(x|\lambda\right)\sim\left(\left(\frac{e^{t}-\lambda}{1-\lambda}\right)^{r},t\right).\label{eq:52}
\end{equation}

Let us assume that
\begin{equation}
HB_{n}^{\left(\nu,k\right)}\left(x\right)=\sum_{m=0}^{n}C_{n,m}H_{m}^{\left(r\right)}\left(x|\lambda\right).\label{eq:53}
\end{equation}

By (\ref{eq:21}), we get
\begin{align}
C_{n,m} & =\frac{1}{m!}\left\langle \left.\left(\frac{e^{t}-\lambda}{1-\lambda}\right)^{r}t^{m}\right|e^{-\frac{\nu t^{2}}{2}}\frac{\li[k]\left(1-e^{-t}\right)}{1-e^{-t}}x^{n}\right\rangle \label{eq:54} \\
 & =\frac{\left(n\right)_{m}}{m!\left(1-\lambda\right)^{r}}\left\langle \left.\sum_{l=0}^{r}\dbinom{r}{l}\left(-\lambda\right)^{r-l}e^{lt}\right|HB_{n-m}^{\left(\nu,k\right)}\left(x\right)\right\rangle  \nonumber \\
 & =\dbinom{n}{m}\frac{1}{\left(1-\lambda\right)^{r}}\sum_{l=0}^{r}\dbinom{r}{l}\left(-\lambda\right)^{r-l}\left\langle \left.1\right|e^{lt}HB_{n-m}^{\left(\nu,k\right)}\left(x\right)\right\rangle  \nonumber\\
 & =\frac{\tbinom{n}{m}}{\left(1-\lambda\right)^{r}}\sum_{l=0}^{r}\dbinom{r}{l}\left(-\lambda\right)^{r-l}HB_{n-m}^{\left(\nu,k\right)}\left(l\right)  \nonumber
\end{align}

Therefore, by (\ref{eq:53}) and (\ref{eq:54}), we obtain the following
theorem.

\begin{thm}

For $n,r\in\mathbb{Z}_{\ge0}$, we have
\[
HB_{n}^{\left(\nu,k\right)}\left(x\right)=\frac{1}{\left(1-\lambda\right)^{r}}\sum_{m=0}^{n}\dbinom{n}{m}\left\{ \sum_{l=0}^{r}\dbinom{r}{l}\left(-\lambda\right)^{r-l}HB_{n-m}^{\left(\nu,k\right)}\left(l\right)\right\} H_{m}^{\left(r\right)}\left(x|\lambda\right).
\]

\end{thm}


\bigskip
ACKNOWLEDGEMENTS. This work was supported by the National Research Foundation of Korea(NRF) grant funded by the Korea government(MOE)\\
(No.2012R1A1A2003786 ).
\bigskip

$\,$

\noindent Department of Mathematics, Sogang University, Seoul 121-742,
Republic of Korea

\noindent e-mail:dskim@sogang.ac.kr\\

\noindent Department of Mathematics, Kwangwoon University, Seoul 139-701,
Republic of Korea

\noindent e-mail:tkkim@kw.ac.kr

\begin{thebibliography}{10}
\bibitem{key-1} S. Araci, M. Acikgoz, \emph{A note on the Frobenius-Euler
numbers and polynomials associated with Bernstein polynomials},\emph{
}Adv. Stud. Contemp. Math. 22 (2012), no.3, 399-406.

\bibitem{key-2} R. Dere, Y. Simsek, \emph{Application of umbral algebra
to some special polynomials}, Adv. Stud. Contemp. Math. 22 (2012),
no. 3, 433-438.

\bibitem{key-3} D. Ding, J. Yang, \emph{Some identities related to
the Apostol-Euler and Apostol-Bernoulli polynomials}, Adv. Stud. Contemp.
Math. 20 (2010), no. 1, 7-21.

\bibitem{key-4} H. Ozden, I. N. Cangul, Y. Simsek, \emph{Remarks
on $q$-Bernoulli numbers associated with Daehee numbers}, Adv. Stud.
Contemp. Math. 18 (2009), no. 1, 41-48.

\bibitem{key-5} M. Kaneko, \emph{Poly-Bernoulli numbers}, J. Th'eor.
Nombres Bordeaux 9(1997), no. 1, 221-228.

\bibitem{key-6} D. S. Kim, T. Kim, D. V. Dolgy, S. H. Rim, \emph{Some
new identities of Bernoulli, Euler and Hermite polynomials arising
from umbral calculus}, Adv. Difference Equ. 2013 (2013), 2013:73.

\bibitem{key-7} D. S. Kim, T. Kim, \emph{Some identities of Frobenius-Euler
polynomials arising from umbral calculus}, Adv. Difference Equ. 2012
(2012), \#196.

\bibitem{key-8} D. S. Kim, T. Kim, S. H. Lee, \emph{A note on poly-Bernoulli
polynomials arising from umbral calculus}, Adv. Studies Theor. Phys.,
7 (2013), no. 15, 731-744.

\bibitem{key-9} D. S. Kim, T. Kim, S. H. Lee, \emph{Poly-Cauchy numbers
and polynomials with umbral calculus viewpoint}, Int. Journal of Math.
Analysis, 7 (2013), 2235-2253.

\bibitem{key-10} D. S. Kim, T. Kim, S. H. Lee, \emph{Higher-order
Cauchy of the first kind and Poly-Cauchy of the first kind mixed type
polynomials}, Adv. Stud. Contemp. Math. 23 (2013), 543-554.

\bibitem{key-11} D. S. Kim, T. Kim, \emph{Some identities of Bernoulli
and Euler polynomials arising from umbral calculus}, Adv. Stud. Contemp.
Math. 23 (2013), no. 1, 159-171.

\bibitem{key-12} B. Kurt, Y. Simsek, \emph{On Hermite based Genocchi
polynomials}, Adv. Stud. Contemp. Math. 23 (2013), no. 1, 13-17.

\bibitem{key-13} S. Roman, \emph{The umbral calculus}, Pure and Applied
Mathematics, 111,Academic Press, Inc. {[}Harcourt Brace Jovanovich,
Publishers{]}, New York, 1984, $x$+193 pp. ISBN: 0-12-594380-6.

\bibitem{key-14} S. Roman, G.-C. Rota, \emph{The umbral calculus,
}Advances in Math. 27 (1978), no. 2, 95-188.

\bibitem{key-15} S. H. Rim, J. Jeong, \emph{On the modified $q$-Euler
numbers of higher order with weight, }Adv. Stud. Contemp. Math. 22
(2012), no. 1, 93-98.

\bibitem{key-16} Y. Simsek, \emph{Generating functions of the twisted
Bernoulli numbers and polynomials associated with their interpolation
functions}, Adv. Stud. Contemp. Math. 16 (2008), no. 2, 251-278.

\end{thebibliography}
\end{document}